\documentclass[10pt, a4paper, reqno]{amsart}

\usepackage{amscd}
\usepackage{dsfont}

\usepackage{booktabs}
\usepackage{changepage}

\usepackage{amsmath}
\usepackage{amstext}
\usepackage{amsthm}
\usepackage{amssymb}
\usepackage{amsopn}
\usepackage{amssymb}
\usepackage{amsxtra}
\usepackage{amsfonts}
\usepackage{fancyhdr, xcolor}
\usepackage{paralist, epsfig}
\usepackage[mathcal]{euscript}
\usepackage[english]{babel}
\usepackage[latin1]{inputenc}

\makeatletter
\g@addto@macro\normalsize{%
  \setlength\abovedisplayskip{10pt}
  \setlength\belowdisplayskip{10pt}
  \setlength\abovedisplayshortskip{5pt}
  \setlength\belowdisplayshortskip{8pt}
}

\usepackage{pgfplots}
\pgfplotsset{compat=1.3}

\newcommand\clipright[1][white]{
  \fill[#1](current axis.south east)rectangle(current axis.north-|current axis.outer east);
  \pgfresetboundingbox
  \useasboundingbox(current axis.outer south west)rectangle([xshift=.5ex]current axis.outer north-|current axis.east);
}

\definecolor{mycolor}{rgb}{0.02,0.4,0.7}

\usepackage{setspace}
\usepackage{color}
\newcommand{\red}[1]{\textcolor{black}{#1}}
\newcommand{\T}{\operatorname{T}}
\newcommand{\F}{\operatorname{F}}

\definecolor{grey}{gray}{.3}

\newcommand{\Bigsum}[2]{\ensuremath{\mathop{\textstyle\sum}_{#1}^{#2}}}
\begin{document}

\title{An outlook on self-assessment of homework\\assignments in higher mathematics education}

\author{Sarah Beumann\hspace{0.5pt}\MakeLowercase{$^{\text{1}}$} and Sven-Ake Wegner\hspace{0.5pt}\MakeLowercase{$^{\text{2}}$}}

\renewcommand{\thefootnote}{}
\hspace{-1000pt}\footnote{\hspace{5.5pt}2010 \textit{Mathematics Subject Classification}: Primary 97D60; Secondary 97D40, 97B40\vspace{1.6pt}}

\hspace{-1000pt}\footnote{\hspace{5.5pt}\textit{Key words and phrases}: self-assessment, self-evaluation, self-regulation, mathematics education, tertiary education. \vspace{1.6pt}}




\hspace{-1000pt}\footnote{\hspace{0pt}$^{1}$\,University of Wuppertal, School of Mathematics and Natural Sciences, Gau\ss{}stra\ss{}e 20, 42119 Wuppertal, Germany,\linebreak\phantom{x}\hspace{1.2pt}Phone: +49\hspace{1.2pt}(0)\hspace{1.2pt}202\hspace{1.2pt}/\hspace{1.2pt}439\hspace{1.2pt}-\hspace{1.2pt}5192, E-Mail: beumann@uni-wuppertal.de.\vspace{1.6pt}}

\hspace{-1000pt}\footnote{\hspace{0pt}$^{2}$\,Corresponding author: Teesside University, School of Science, Engineering \&{} Design, Stephenson Building, Middles-\linebreak\phantom{x}\hspace{1.2pt}brough, TS1\;3BX, United Kingdom, phone: +44\,(0)\,1642\:434\:82\:00, e-mail: s.wegner@tees.ac.uk.}

\begin{abstract}
We discuss first experiences with a new variant of self-assessment in higher mathematics education. In our setting, the students of the course have to \red{mark} a part of their homework assignments themselves and they receive the corresponding credit without that any later changes are carried out by the teacher. In this way we seek to correct the imbalance between student-centered learning arrangements and assessment concepts that keep the privilege to grade \red{(or mark)} completely with the teacher. \red{We present results in the form of student feedback from a course on functional analysis for 3rd and 4th year students. Moreover we analyze marking results from two courses on real analysis. Here, we compare tasks marked by the teacher and tasks marked by the students.}
\end{abstract}

\maketitle

\vspace{-10pt}


\section{Self-assessment}\label{SEC-INT}

\smallskip

In \red{recent} years the possibilities to access mathematical knowledge \red{have} increased significantly due to the digitalization of classical media like textbooks, exercises or model solutions, and due to concepts \red{such as} blogs, internet forums, online-available video-taped lectures etc. Modern teaching methods aim \red{to facilitate} the latter to improve students' learning success. \red{They} achieve this by using student-centered learning arrangements such as problem-based learning, research-based learning or other methods that give the students more freedom, but also assign more responsibility to them for their own learning outcome. However, when it comes to an \red{assessment}, often classical instruments, like graded homework assignments, weekly quizzes or closed-book exams, prevail. The philosophy behind this paper is the idea \red{of improving} the imbalance between learning arrangements and assessment by sharing, to some extent, the \red{teachers'} privilege to grade \red{(or to mark)} with the students. Our concrete aim is to strengthen the students' sense \red{of} being responsible for their own learning process by sharing with them the control. This in turn encourages the students to employ the advantages of digitalization to increase their own learning success. In particular, \red{they} no longer feel the need to hide the sources of their ideas from the teacher, but can themselves evaluate their personal gain in knowledge, skills and competencies that they \red{have} extracted from these sources. \red{The latter is} a very important aspect of modern student-centered education.

\medskip

The idea \red{of sharing} the control over the learning process with the students is neither new nor a concept that can easily be realized in the classroom. Indeed, Klenovski's \cite[p.~161]{K1995} quotation \red{from a 1994 interview with a college teacher has lost nothing of its relevance:}

\medskip

\begin{adjustwidth}{1.5cm}{1.5cm}\textquotedblleft{}Students have to learn that it's their course, their learning and they have to take some control\dots{}it's hard for some students because they want you to take control.\textquotedblright{}
\end{adjustwidth}

\medskip

 However, from the mid 90s on, different realizations of the idea have been surveyed in many areas of education such as chemistry \cite{Davey2015, K1995}, mathematics and statistics \cite{Handbook, OS2004, RHGR2001, RHGR2002, BAZF2004, ST1996}, music \cite{H2010}, narrative writing \cite{RRH1999}, and with students of different \red{ages} and school \red{types} such as elementary school \cite{BAZF2004}, middle school \cite{H2010, RRH1999}, high school and college \cite{ST1996}, to list only a sample. Some of these surveys mention a positive impact on the students' achievement \cite{FF1994, RHGR2002, BAZF2004, SC2005}\red{;} some mention no impact \cite{H2010, RHGR2001}\red{;} some point out that self-assessment is not always precise \cite{BBHC2012, Davey2015}. A positive influence on meta-competencies like self-efficacy \cite{RHGR2002}, self-confidence \cite{OS2004}, active learning and motivation \cite{FF1994} as well as critical thinking and the ability to \red{reflect on own work} \cite{Cooper2006} is \red{mentioned}. In \cite{Handbook} it is pointed out that appropriate beliefs about mathematics and \red{mathematical} learning are an important precondition.

\medskip

In the papers cited above rather different approaches are outlined about how to share control with the students in a concrete classroom situation. In this paper we follow mostly the ideas of Klenovski \cite{K1995} who used the two notions of \textit{self-evaluation} and \textit{self-assessment}. Indeed, Klenovski \cite[p.~155--160]{K1995} identifies \textquotedblleft{}three key dimensions of the student self-evaluation process~[\dots\!]:~the use of criteria by students to self-evaluate their own learning~[\dots\!], the interactive dialogue~[\dots\!]~between student and teacher, during the ana\-ly\-sis of the student's self-evaluation, [and] the ascription of a grade by the students for their own work.\textquotedblright{} Klenovski \cite[p.~147]{K1995} states that \textquotedblleft{}self-evaluation~[\dots\!]~is broader than self-assessment in that the student is engaged in more than just deciding what grade he or she should get.\textquotedblright{} It appears to us that in the classroom situations surveyed by Klenovski the students did not have the final authority about the grade, but that the teacher could intervene \cite[second interview on p.~159]{K1995}, or an intervention by peer-learners was possible \cite[interview on p.~158]{K1995}. In our experiments it is essential that the students ascribe their own grades \red{(or marks)} \textit{without \red{the} intervention} of a second party. For this reason we stick below to the word \textit{self-assessment} although, of course, the use of criteria, and a dialogue about assessment, are important in our setting as well. Our \red{incentive} behind this concept of self-assessment---which differs \red{from} our knowledge \red{of} all concepts discussed so far in the literature---is the following.

\medskip

\begin{compactitem}

\item[1.] Self-assessment allows \red{us} to give meta-tasks to the students that cannot be \red{marked} by the teacher. \red{Examples could be} to repeat some topic from the last \red{year's} course or to practise a method \textquotedblleft{}until the students master it\textquotedblright{}.

\smallskip

\item[2.] Self-assessment allows \red{us} to give extra tasks to the students, and to grant credit for working on these tasks, without \red{the school having to pay staff that carries out the marking}.

\smallskip

\item[3.] \red{Self-assessment helps to illustrate} that checking the validity of a proof is not a formal and fail-safe procedure but requires careful work and \red{may depend} on personal taste. \red{This is for example the case} when it comes to the amount of details that are given and the strategy that is pursued. \red{In this sense self-assessment generates appropriate beliefs about mathematics.}

\smallskip

\item[4.] Self-assessment transfers to the students\red{, for a moment,} the full \red{responsibility for} their grading \red{(or marking)} and thus \red{fosters} the development of the earlier mentioned meta-competencies\red{---}like self-efficacy, self-confidence, and motivation\red{---}compared to \red{situations} in which students participate in the evaluation but the final grading \red{(or marking)} is done by the teacher.

\smallskip

\item[5.] Self-assessment encourages the \red{students not just to maximize the teacher-assigned grade but  to learn mathematics on a level of deep understanding.}  

\end{compactitem}

\medskip

Let us give two examples of authentic classroom situations that illustrate our \red{incentive behind this article}. In situation 1 a student \red{kept} asking for help with an exercise until the teacher solved the whole task for the student. As the solution is now of course correct, the teacher \red{assigned}, after it \red{was} handed in, the \red{maximum} number of \red{marks}. The student's learning progress might however be poor as mathematics is not about applying internalized techniques to well-known problems, but about finding new techniques to solve unknown problems---which students \red{only} learn by solving problems on their own. In situation 2 the student hands in a solution copied from a book or from the internet. From the solution the teacher can see that it was copied without any understanding, e.g., as it follows a naming convention different from the lecture, or as the notation is completely different from that on the problem sheet. As the math is however correct, the teacher feels that he cannot deduct much from the full score. The student's learning progress is\red{,} however\red{,} more or less zero. Our initial idea was that giving the power and duty of \red{marking} to the students \red{in such situations} could result in a change of their beliefs. It could help the students to reconsider their strategies and become aware of their own responsibility---for their learning progress and for the mathematical work that they produce.

\medskip

Let us mention that our basic idea \red{of giving} more control to the students \red{in order to} improve the learning process, is also the leitmotif in Klenovski's paper \cite{K1995}. His findings \cite[p.~161f]{K1995} support the latter statement but also point out that pedagogical change is needed and implementations of the concept have to be \red{further studied}. The first results explained below confirm that our new concept of intervention-free self-assessment can be applied successfully in higher mathematics education. On the other hand they also identify drawbacks and obstructions. This paper is intended as a small preview and an invitation to other university teachers to contribute with their ideas and experience to the development of self-assessment in mathematics.

\medskip
\section{A pilot study---first results on self-assessment}\label{SEC-EXP}

In this section we outline first experiences with our concept of self-assessment by presenting students' feedback and \red{the marks} of two homework assignments. \red{We} compare \red{the} results of \red{parts that were} assessed by the teacher with \red{parts} that \red{were} assessed by the students. \red{We} present and discuss some selected feedback that gives \red{insight into} the students' beliefs about their role in the learning and assessment process.

\bigskip

\begin{center}
{\sc2.1.~Homework assignments in higher mathematics}
\end{center}

\medskip

The first experience \red{of the authors} with self-assessment was the spontaneous idea to assign the review of topics that had been covered in a previous course as \red{a homework assignment}. \red{In order to underline that we wanted this to be understood as a serious task we decided to put it in the following form as one of four tasks on the weekly exercise sheet.}

\medskip

\begin{adjustwidth}{1.5cm}{1.5cm}
\textbf{Exercise 1.~(5 \red{marks})} Review the construction of the Lebesgue integral, the dominated convergence theorem and the monotone convergence theorem. Maybe it is helpful to browse the appendix of the book \cite{Werner} by D.~Werner.
\end{adjustwidth}

\medskip

This task was given in the middle of a 14-\red{week} course on the foundations of functional analysis taught in 2012 with approximately 20 students in their 3rd and 4th year. Each exercise sheet contained four tasks for which solutions had to be handed in \red{and that were usually marked by the teacher}. On this particular sheet only three tasks \red{required a} solution. For the forth one, Exercise 1, the students were required to self-assess their achievement and to indicate the score on the submission. \red{When we handed out the sheet the students appeared very surprised and suspicious because they were not used to exercises of this type.} \red{Many of them} did not award themselves the full amount of \red{five marks}. \red{Indeed,} they assumed that \red{we would carry out some kind of} \textquotedblleft{}double checking\textquotedblright{}, like an oral examination during the recitation, \red{if they assign themselves a high score}. After the semester we received the following feedback by one of the students.

\medskip

\begin{adjustwidth}{1.5cm}{1.5cm}
\textquotedblleft{}The exercise to recall the introduction (definition and main properties) of the Lebesgue integral and to give yourself \red{marks} on the basis of your comprehension is meaningful and helpful as well. First, one recalls the content carefully which leads to a deep understanding, and second the already rehearsed content anchors in memory. Since one gives \red{marks} on the basis of comprehension, you repeat the content carefully to \textquoteleft{}obtain\textquoteright{} a good score. Indeed, in order to avoid an embarrassing situation where the tutor checks that the number of \red{marks} is inappropriate, you think twice of how many points are eligible.\textquotedblright{}
\end{adjustwidth}

\medskip

\red{We mention that Exercise 1, as stated above, was the only self-assessed assignment in this course. The five marks correspond to approximately 2\% of the total score of 260 marks that the students could achieve on the 13 exercise sheets.}

\medskip

Our second experience with self-assessment was the following. \red{During a 14-week course on real analysis, taught in 2013 for 1st year students, we gave the following two exercises.} \red{Both} were given as additional exercises and were credited with 20 \red{marks}. \red{The} total of regular \red{marks} was 480. The self-assessment homework thus \red{counted} as approximately 4\% extra credit.

\medskip

\begin{adjustwidth}{1.5cm}{1.5cm}
\textbf{Exercise 2.~(10 \red{marks})} Become confident with handling sequences and with computing their limits, e.g., by working on the exercises from the additional worksheet on the course's website.
\end{adjustwidth}

\medskip

The additional worksheet contained 46 \red{sequences for which the limits had to be computed}. The second exercise refers to the following theorem that establishes some basic rules for computations with convergent series.

\medskip

\begin{adjustwidth}{1.5cm}{1.5cm}
\textbf{Theorem 1.} Let $(a_k)_{k\geqslant0}$, $(b_k)_{k\geqslant0}\subseteq\mathbb{R}$ and $\lambda\in\mathbb{R}$ be given.

\smallskip

\begin{compactitem}

\item[(i)] We have
$$
\displaystyle\Bigsum{k=0}{\infty}a_k+\lambda\Bigsum{k=0}{\infty}b_k = \Bigsum{k=0}{\infty}a_k+\lambda b_k
$$
provided that the two series on the left are convergent.

\smallskip

\item[(ii)] Assume that there exists $k_0\geqslant0$ such that $a_k=b_k$ holds for all $k\geqslant k_0$. Then the series over all $a_k$'s converges if and only if the series over all $b_k$'s converges.

\smallskip

\item[(iii)] Let the series over the $a_k$'s be convergent and let $(j_k)_{k\geqslant0}$ with $j_k\nearrow\infty$ and $j_0=-1$ be given. Then the following series
$$
\displaystyle\Bigsum{k=0}{\infty}a_{j_k+1}+\cdots+a_{j_{k+1}}\vspace{-5pt}
$$
is also convergent. The converse is \red{false}.
\end{compactitem}
\end{adjustwidth}

\medskip

During the \red{lectures we presented Theorem 1} without its proof. The exercise was then as follows.

\medskip

\begin{adjustwidth}{1.5cm}{1.5cm}
\textbf{Exercise 3.~(10 \red{marks})} Make sure that you are able to prove the rules for computations with convergent series given in Theorem 1, e.g., by \red{giving} all or a suitable selection of the proofs yourself.
\end{adjustwidth}

\medskip

As in \red{Exercise 1}, we asked the students to award themselves the corresponding \red{marks} and to \red{indicate} the score on their \red{submissions}. \red{They} did neither get a model solution or a \red{marking} scheme. \red{This reflects one of the main incentives for self-assessment mentioned in Section \ref{SEC-INT}: Leaving the proofs completely to the students will grow their ability to evaluate if a mathematical argument is correct or not \textit{by themselves}.}

\medskip

\red{We mentioned that} Exercise 2 appeared as an additional task on one of the homework sheets. \red{On this sheet, four exercises were graded by the teacher and one exercise was subject to self-assessment}. The following table shows \red{the averages of the teacher-assessed part and the averages of the self-assessed part}. It is eye-catching that in this case the average of the teacher assessment is \red{approximately} 53\% whereas the average of the self-assessment is \red{approximately} 75\%.

\medskip

\begin{center}
\begin{minipage}[h]{441pt}
\begin{center}
\begin{tabular}{lccc}
    \toprule
         & Sample Size$\phantom{\displaystyle X^{X^X}}$ & \hspace{10pt}Mean\hspace{10pt} & \hspace{22pt}Standard Deviation\vspace{2pt}\\ \midrule\vspace{3pt}

Teacher's assessment$\phantom{\displaystyle X^{X^X}}$& 45$\phantom{\displaystyle X^{X^X}}$ &  21.09&\hspace{22pt}6.73 \\
Self-assessment\hspace{0pt} & 45$\phantom{\displaystyle X^{X^X}}$ & 7.31 &  \hspace{22pt}3.04 \\\bottomrule
  \end{tabular}

\bigskip

{

\small

\textbf{Table 1.}~min=0, max=40 for exercises graded by the teacher,\\ min=0, max=10 for the exercise graded by the students.

}
\end{center}
\end{minipage}
\end{center}

\medskip

The distribution of the teacher-assessed tasks looks Gaussian-like if one ignores the 13\% of students that obtained less than or equal to 10 out of 40 \red{marks}. \red{In this course 50\% of the marks on the sheets were sufficient to be admitted to the final exam. The grade for the course depended only on this exam. In view of this the latter seems reasonable and expectable.}

\begin{center}
\begin{minipage}[h]{441pt}
\begin{center}
\begin{minipage}[h]{300pt}
\hspace{14pt}\begin{tikzpicture}
    \begin{axis}[
        ymajorgrids,
        xmajorgrids,
        grid style={white,thick},
        axis on top,
        /tikz/ybar interval,
        tick align=outside,
        ymin=0,
        axis line style={draw opacity=0},
        tick style={draw=none},
        enlarge x limits=false,
        height=5cm,
        title={},
        title style={font=\Large},
        xlabel={},
        ylabel={},
        ytick={0,1,2,3,4,5,6,7,8,9,10,11,12,13},
        scaled ticks=false,
        yticklabels={},
        xticklabels={$0$, $4$, $8$, $12$, $16$, $20$,
                           $24$, $28$, $32$, $36$, $40$},
        width=\textwidth,
        xtick=data,
        label style={font=\large},
        ticklabel style={
            inner sep=1pt,
            font=\small
        },
        nodes near coords,
        every node near coord/.append style={
            fill=white,
            anchor=center,    
            shift={(13pt,4pt)},
            inner sep=0,
            font=\footnotesize,
            rotate=0},
            ]
    \addplot[mycolor!80!white, fill=mycolor, draw=none] coordinates {(0,  0) (1,  0) (2,  7) (3,  3) (4,  9) (5,  13) (6,  8) (7,  3) (8,  2) (9,  0) (10,  0)};
    \end{axis}
    \clipright
        \draw (-0.1,0) -- (9.15,0);
\end{tikzpicture}
\end{minipage}

\bigskip

{

\small

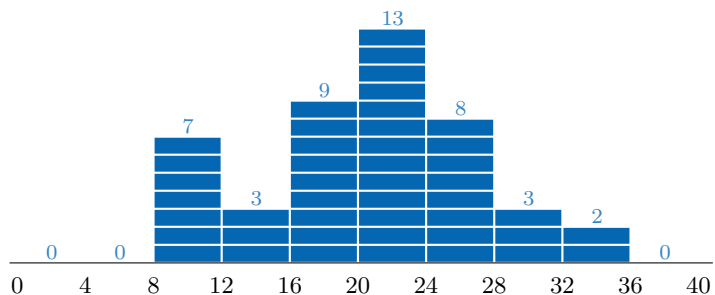
\textbf{Figure 1.} Teacher's assessment: 0 students were awarded $n\in[0,8]$\\points, 7 students were awarded $n\in(8,12]$  points, etc.

}
\end{center}
\end{minipage}
\end{center}

\medskip

The distribution of the student-assessed tasks looks completely different and has a higher average. \red{We prefer} to be careful with drawing conclusions, since we compare exercises on different topics and with different levels of difficulty. \red{It is, however, again eye-catching} that 53\% of the students awarded themselves the full 10 \red{marks}, whereas 18\% awarded themselves zero \red{marks}.

\begin{center}
\begin{minipage}[h]{441pt}
\begin{center}
\begin{minipage}[h]{300pt}
\hspace{15pt}\begin{tikzpicture}
    \begin{axis}[
        ymajorgrids,
        xmajorgrids,
        grid style={white,thick},
        axis on top,
        /tikz/ybar interval,
        tick align=outside,
        ymin=0,
        axis line style={draw opacity=0},
        tick style={draw=none},
        enlarge x limits=false,
        height=4cm,
        title={},
        title style={font=\Large},
        xlabel={},
        ylabel={},
        ytick={0,2,4,6,8,10,12,14,16,18,20,22,24},
        scaled ticks=false,
        yticklabels={},
        xticklabels={$0$, $1$, $2$, $3$, $4$, $5$,
                           $6$, $7$, $8$, $9$, $10$},
        width=\textwidth,
        xtick=data,
        label style={font=\large},
        ticklabel style={
            inner sep=1pt,
            font=\small
        },
        nodes near coords,
        every node near coord/.append style={
            fill=white,
            anchor=center,    
            shift={(13pt,4pt)},
            inner sep=0,
            font=\footnotesize,
            rotate=0},
            ]
    \addplot[mycolor!80!white, fill=mycolor, draw=none] coordinates {(0,  8) (1,  0) (2,  0) (3,  0) (4,  1) (5,  5) (6,  3) (7,  3) (8,  1) (9,  24) (10,  0)};
    \end{axis}
    \clipright
        \draw (-0.1,0) -- (9.15,0);
\end{tikzpicture}
\end{minipage}

\bigskip

{

\small

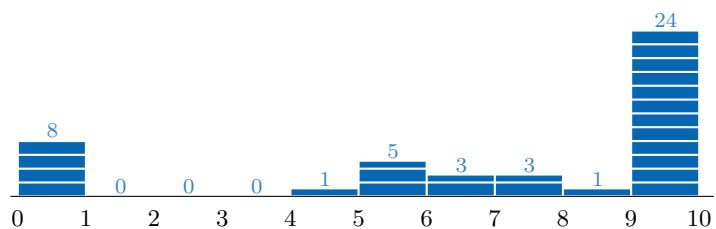
\textbf{Figure 2.} Self-assessment: 8 students awarded themselves $n\in[0,1]$\\points, 0 students awarded themselves $n\in(1,2]$ points, etc.

}

\end{center}
\end{minipage}
\end{center}

\medskip

\red{It seems very interesting and important to us} that among those eight students that assigned themselves \red{zero marks}, only one received 10 of 40 \red{marks} from the teacher. \red{The} other seven received between 17 and 26 \red{out} of 40 and thus \red{scored} around the average value. \red{Among} the 24 students that gave themselves the full 10 \red{marks}, we find five out of those six students that received less than or equal to 10 \red{marks} in the \red{teacher-assessed part}. This suggests that \red{weak students in particular} did not assess themselves very honestly. For a further development of self-assessment techniques this effect has to be taken into account. More experiments \red{are needed to see} if \red{the latter} is a general trend or if the students in the long \red{term} will \red{assess} themselves \red{in a reasonable fashion}.

\medskip

\red{The third experiment on self-assessment was part of a 12-week course on real analysis for 1st year students taught in 2018. We mention that we had a very small group of only seven students and thus an atmosphere in which the students know each other well and talk much about math, homework, exams etc. The assessment consisted of a final exam and one longer homework assignment in the middle of the course. Both components contributed 50\% to the final grade. The homework assignment consisted of 10 questions. It covered elementary logic, sets, mappings and mathematical induction. One of the 10 questions was the following.
\medskip
\begin{adjustwidth}{1.5cm}{1.5cm}
\textbf{Exercise 4.~(10 marks)} Become confident with using truth tables by verifying a suitable sample the following statements:\vspace{2pt}
\begin{compactitem}
\item[\phantom{00}(1)] $A\wedge \T \Leftrightarrow A$, $A\vee \F \Leftrightarrow A$\vspace{2pt}
\item[\phantom{00}(2)] $A\vee \T \Leftrightarrow T$, $A\wedge \F \Leftrightarrow F$\vspace{2pt}
\item[\phantom{00}(3)] $A\vee A \Leftrightarrow A$, $A\wedge A \Leftrightarrow A$\vspace{2pt}
\item[\phantom{00}(4)] $\neg{}(\neg A) \Leftrightarrow A$\vspace{2pt}
\item[\phantom{00}(5)] $A\vee B \Leftrightarrow B\vee A$, $A\wedge B \Leftrightarrow B\wedge A$\vspace{2pt}
\item[\phantom{00}(6)] $A\vee (B\vee C) \Leftrightarrow (A\vee B)\vee C$, $A\wedge (B\wedge C) \Leftrightarrow (A\wedge B)\wedge C$\vspace{2pt}
\item[\phantom{00}(7)] $A\vee (B\wedge C) \Leftrightarrow (A\vee B)\wedge (A\vee C)$, $A\wedge (B\vee C) \Leftrightarrow (A\wedge B)\vee (A\wedge C)$ \vspace{2pt}
\item[\phantom{00}(8)] $\neg(A\wedge B) \Leftrightarrow \neg A \vee \neg B$, $\neg(A\vee B) \Leftrightarrow \neg A \wedge \neg B$ \vspace{2pt}
\item[\phantom{00}(9)] $\neg(A\wedge B) \Leftrightarrow \neg A \vee \neg B$, $\neg(A\vee B) \Leftrightarrow \neg A \wedge \neg B$ \vspace{2pt}
\item[\phantom{0}(10)] $(A\Rightarrow B)\Leftrightarrow (\neg A\vee B)$\vspace{2pt}
\item[\phantom{0}(11)] $A\vee \neg A$, $\neg(A\wedge\neg A)$ \vspace{2pt}
\item[\phantom{0}(12)] $[(A\Rightarrow B)\wedge \neg B]\Rightarrow \neg A$ \vspace{2pt}
\item[\phantom{0}(13)] $[(A\Rightarrow B)\wedge (B\Rightarrow C)]\Rightarrow (A\Rightarrow C)$ \vspace{2pt}
\item[\phantom{0}(14)] $(A\wedge B)\Rightarrow A$, $(A\wedge B)\Rightarrow B$  \vspace{2pt}
\item[\phantom{0}(15)] $A\Rightarrow (A\vee B)$, $B\Rightarrow (A\vee B)$\vspace{2pt}
\item[\phantom{0}(16)] $(A\Leftrightarrow B)\Leftrightarrow [(A\Rightarrow B)\wedge(B\Rightarrow A)]$ \vspace{2pt}
\item[\phantom{0}(17)] $(A\Rightarrow B)\Leftrightarrow(\neg B\Rightarrow\neg A)$\vspace{2pt}
\item[\phantom{0}(18)] $[(A\vee B)\wedge\neg A]\Rightarrow B$ \vspace{2pt}
\item[\phantom{0}(19)]  $[(\neg A\wedge B)\Rightarrow F]\Rightarrow (A\Rightarrow B)$ \vspace{2pt}
\item[\phantom{0}(20)] $[(A\Rightarrow B)\wedge A]\Rightarrow B$\vspace{3pt}
\end{compactitem}
Indicate the number of marks on your submission. Don't hand in any truth table!
\end{adjustwidth}
\medskip
Exercise 4 contributed 5\% to the final mark. The design was similar to Exercise 2, where we gave 46 sequences to practise the computation of limits. However, we point out that the computation of these limits in most cases involved a certain trick, like applying an estimate, or combining two previous limits in a suitable way. In contrast to this, Exercise 4 was much more straightforward and can be completed---once the principle is understood---by a rather \textquotedblleft{}mechanical procedure\textquotedblright{}.}

\medskip

\red{In the table below we compare again the grading results of the self-assessed part with the teacher-assessed part.}

\medskip

\begin{center}
\begin{minipage}[h]{441pt}
\begin{center}
\red{\begin{tabular}{lccc}
    \toprule
         & Sample Size$\phantom{\displaystyle X^{X^X}}$ & \hspace{10pt}Mean\hspace{10pt} & \hspace{22pt}Standard Deviation\vspace{2pt}\\ \midrule\vspace{3pt}
Teacher's assessment$\phantom{\displaystyle X^{X^X}}$& 7$\phantom{\displaystyle X^{X^X}}$ & 56.86 &\hspace{22pt}16.00 \\
Self-assessment\hspace{0pt} & 7$\phantom{\displaystyle X^{X^X}}$ & 7.86 &  \hspace{22pt}3.67 \\\bottomrule
  \end{tabular}}

\bigskip

{

\small

\red{\textbf{Table 2.}~min=0, max=90 for exercises graded by the teacher,\\ min=0, max=10 for the exercise graded by the students.}

}
\end{center}
\end{minipage}
\end{center}

\medskip

\red{In our small group of seven students the average of the tasks assessed by the teacher was with 63\% lower than the 79\% of the self-assessed part. This was also the case with Exercise 2. The correlation between the marks that the students gave themselves and the marks that the teacher gave to them was 0.77 in the current experiment. In the previous experiment the correlation was only 0.05. One might conclude from this that the students' evaluation of their own abilities in this case was closer to the teacher's evalutation of the latter. However, we would like to be cautious here in view of the small group size and  the different types of questions in Exercise 2 and Exercise 4. On the other hand, we are indeed convinced that this last experiment with self-assessment was more successful than the previous one. We recognized that some of the students put much effort into Exercise 4 and indeed did all 31 truth tables. By doing this, they gained not only the desired proficiency with the method. At the same time they gained confidence in their own abilities and handed in their solutions with the good feeling that they really deserve the 10/10 marks that they ascribed to themselves. With a classical design (one or two of the statements listed in Exercise 4 to be handed in and to be marked by the teacher) we could not have achieved this.}

\medskip

\begin{center}
{\sc2.2.~Student's impressions about individual responsibility}
\end{center}

The last experience that we want to discuss here \red{did not involve self-assessment} in the sense of Section \ref{SEC-INT}. \red{It was, however, similar in the sense that} the responsibility to work on homework assignments was \red{completely due to the students}. \red{In contrast to the situations explained above} the \red{marking} was waived \red{completely}. \red{In a} 3rd year course with approximately 10 students and \red{in a} 2nd year course with approximately 50 students we strongly recommended \red{intensive work} on the weekly assignments. \red{We} emphasized that the final exam will be very \red{similar} to the tasks \red{in these assignments}. \red{In} the small course \red{we asked} the students to present their solutions during the exercise sessions. In the large course the solutions were presented by the teacher and later uploaded to the website of the course, as there were too many participants for individual presentations. The grade for both courses was given on the basis of the final exam. During the \red{term} we received much negative feedback. \red{Indeed, most of the other teachers employed homework assessment, quizzes, midterm exams and strict attendance requirements to control the students' engagement}. In view of the exam outcome one can say that our concept completely failed in this context. In the middle of the course we already recognized that only less than one quarter of the students downloaded the exercise sheets before the lesson. The whole situation is very well \red{summarized} by the following feedback \red{comment}.

\medskip

\begin{adjustwidth}{1.5cm}{1.5cm}
\textquotedblleft{}100\% final is\:\dots\,strange\:\dots\,it has good and bad sides. Bad thing is that the students sometimes \textquoteleft{}forget\textquoteright{} about this course for the whole semester, which affects their final preparation.\textquotedblright{}
\end{adjustwidth}

\medskip

\red{From this} one can deduce that the students were indeed aware that they did not assume responsibility for their own learning progress. However, \red{it was us} who did not manage to initiate a change of their learning behavior in \red{this course}. On the other hand, we received \red{the following} positive comment.

\medskip

\begin{adjustwidth}{1.5cm}{1.5cm}
\textquotedblleft{}Learning the subject WITHOUT WORRYING that you fail quiz or midterm and don't have chance to pass the course. Learning with our own pace. Mock exams and homeworks help much. It seems risky and stressful at the end. But I think having too much midterms and quizzes give constant stress which makes hard student life for low-pace studiers.\textquotedblright{}
\end{adjustwidth}

\medskip

\red{This comment} suggests that a paradigm shift \red{might have been} possible, but would \red{had required} a different methodology. \red{Self-assessment---that we unfortunately did not use in this case---could have improved the situation.}


\section{Discussion and Outlook}\label{SEC-DIS}

\smallskip

Self-assessment in \red{the sense of this article} can be used successfully in higher mathematics education. The feedback from our 3rd year course on functional analysis indicated that students assessed themselves honestly or even too cautiously. In the \red{first experiment with 1st year students} the data \red{indicates} that on average students overrated themselves within the self-assessed tasks and that in particular the very weak students did this excessively. Of course it is also possible that the teacher underrated certain students in the non self-assessed tasks. Indeed, it is a key problem \red{of} assessment that \red{the latter is} always subjective and individual. In view of the low weight \red{($\leqslant5\%$)} of the self-assessment tasks, we consider the overrating as a tolerable side-effect. \red{The setting of a small group and a task such as Exercise 4---in which everybody can achieve the full score by hard work---turned out to be very suitable for self-assessment. This setting in particular seems to grow the weaker students' confidence in their own abilities.} \red{We point out that} our concept differs substantially from previous implementations of self-evaluation due to the fact that students actually \red{mark} their own work without \red{interventions} of peers or the teacher. In particular the first and the last comment presented in Section \ref{SEC-EXP} suggest that this \red{amplifies} the belief that an effective learning process has to be designed by teachers and students together. 

\medskip

\red{Our first explorative results also identify drawbacks and obstructions.} The first comment in Section \ref{SEC-EXP}.2 illustrates that it can be \red{very difficult} to achieve that students develop a \red{sense of responsibility}. In certain environments it might even be impossible. \red{Our} experiments highlight that \red{we} cannot expect a priori that students will grade themselves honestly. Therefore sophisticated implementations need to be designed in the future. In order to improve our concept we aim to get an in-depth look into the self-evaluation process itself. It would be desirable to \red{obtain} more information on \textit{how} the students actually ascribe the \red{marks}. However, collecting the \red{students'} solutions and assessing their assessment---even if only for research purposes---might \red{already} influence the self-assessment. \red{It seems to us} that there is no easy or standard way to implement self-assessment.

\medskip

To conclude, we like to mention once more that this small preview \red{is} intended as an invitation to other university teachers to contribute with their ideas and experience to the \red{topic} of self-assessment in mathematics. Larger experiments, that will follow the lines sketched above, are under preparation.

\bigskip

\footnotesize

{\sc Acknowledgements.} The authors would like to thank the referees for many helpful and constructive comments that helped to improve this paper significantly. Moreover, the authors would like to thank B.~Farkas (Wuppertal) who taught the course from which we took the Exercises 1 and 2, and who supported the authors with several valuable comments during the preparation of this article. \red{Finally, the authors would like to thank K.~Turner (Teesside) for many useful advices that helped to improve this article significantly.}

\end{document}